\documentclass{article}%
\usepackage{amsmath}
\usepackage{amsfonts}
\usepackage{amssymb}
\usepackage{graphicx}%
\providecommand{\U}[1]{\protect\rule{.1in}{.1in}}
\newtheorem{theorem}{Theorem}

\newtheorem{definition}[theorem]{Definition}

\newtheorem{remark}[theorem]{Remark}

\begin{document}

\title{Approximation properties of the $q$-Bal\'{a}zs-Szabados operators in the case
$q\geq1$}
\author{N. I. Mahmudov\\Department of Mathematics\\Eastern Mediterranean University \\Gazimagusa, TRNC\\Mersin 10, Turkey\\Email: nazim.mahmudov@emu.edu.tr}
\date{}
\maketitle

\begin{abstract}
This paper deals with approximation properties of the newly defined
$q$-generalization of the Bal\'{a}zs-Szabados operators in the case $q\geq1$.
Quantitative estimates of the convergence and Voronovskaja type theorem are
given. In particular, it is proved that the rate of approximation by the
$q$-Bal\'{a}zs-Szabados ($q>1$) is of order $q^{-n}$ versus $1/n$ for the
classical Bal\'{a}zs-Szabados ($q=1$) operators. The results are new even for
the classical case $q=1$.

\end{abstract}

\section{Introduction}

The goal of the paper is to define a $q$-analogue and study approximation
properties for the rational complex Bal\'{a}zs-Szabados operators given by%
\[
R_{n}\left(  f;z\right)  =\dfrac{1}{\left(  1+a_{n}z\right)  ^{n}}\sum
_{k=0}^{n}f\left(  \dfrac{k}{b_{n}}\right)  \left(
\begin{array}
[c]{c}%
n\\
k
\end{array}
\right)  \left(  a_{n}z\right)  ^{k},
\]
where $a_{n}=n^{\beta-1}$, $b_{n}=n^{\beta}$, $0<\beta\leq\frac{2}{3}$,
$n\in\mathbb{N}$ and $z\in\mathbb{C}$, $z\neq-\frac{1}{a_{n}}$.

In the real form rational operators were introduced and studied in Bal\'{a}zs
\cite{balazs1} and Bal\'{a}zs-Szabados \cite{balazs2}. Totik \cite{totik}
settled the saturation properties of $R_{n}\left(  f\right)  $. Further
studies on these operators in the case of real variable can be found in the
paper Abel-Della Vecchia \cite{abel}. They studied the complete asymptotic
expansion for operators $R_{n}\left(  f\right)  $ as $n\rightarrow\infty$.
Approximation properties of the complex Bal\'{a}zs-Szabados operators were
studied in Gal \cite{galbook1}. The $q$-analoque of these operator was given
by Do\u{g}ru who investigated statistical approximation properties of
$q$-Bal\'{a}zs-Szabados operators \cite{dogru}. The approximation properties
of the complex $q$-Bal\'{a}zs-Szabados operators is studied in \cite{ozkan}.

We introduce some notations and definitions of $q$-calculus, see \cite{AAR},
\cite{agarwal2}. Let $q>0.$ For any $n\in\mathbb{N}\cup\left\{  0\right\}  $,
the $q$-integer $\left[  n\right]  _{q}$ is defined by%
\[
\left[  n\right]  _{q}:=\left\{
\begin{tabular}
[c]{ll}%
$\left(  1-q^{n}\right)  /\left(  1-q\right)  ,$ & $q\neq1$\\
& \\
$n,$ & $q=1$%
\end{tabular}
\ \ \right.  ,\ \ \ \left[  0\right]  _{q}:=0;
\]
and the $q$-factorial $\left[  n\right]  _{q}!$ by%
\[
\left[  n\right]  _{q}!:=\left[  1\right]  _{q}\left[  2\right]
_{q}...\left[  n\right]  _{q},\ \ \ \ \left[  0\right]  _{q}!:=1.
\]
For integers $0\leq k\leq n$, the $q$-binomial coefficients are defined by%
\[
\left[
\begin{array}
[c]{c}%
n\\
k
\end{array}
\right]  _{q}:=\frac{\left[  n\right]  _{q}!}{\left[  k\right]  _{q}!\left[
n-k\right]  _{q}!}.
\]
For fixed $q>0,\ q\neq1$, we denote the $q$-derivative $D_{q}f\left(
z\right)  $ of $f$ by%
\[
D_{q}f\left(  z\right)  =\left\{
\begin{tabular}
[c]{lll}%
$\frac{f\left(  qz\right)  -f\left(  z\right)  }{\left(  q-1\right)  z},$ &
$z\neq0,$ & \\
&  & \\
$f^{\prime}\left(  0\right)  ,$ & $z=0.$ &
\end{tabular}
\ \ \ \ \right.
\]

Let us introduce a $q$-Bal\'{a}zs-Szabados operator.

\begin{definition}
Let $q>0$. For $f:\left[  0,\infty\right)  \rightarrow\mathbb{R}$ we define
the Bal\'{a}zs-Szabados operator based on the $q$-integers as follows.%
\begin{equation}
R_{n,q}\left(  f;x\right)  =\dfrac{1}{\left(  1+a_{n}x\right)  ^{n}}\sum
_{k=0}^{n}f\left(  \dfrac{\left[  k\right]  _{q}}{b_{n}}\right)  \left[
\begin{array}
[c]{c}%
n\\
k
\end{array}
\right]  _{q}\left(  a_{n}x\right)  ^{k}%
{\displaystyle\prod\limits_{s=0}^{n-k-1}}
\left(  1+\left(  1-q\right)  \left[  s\right]  _{q}a_{n}x\right)  ,
\label{bs}%
\end{equation}
where $a_{n}=\left[  n\right]  _{q}^{\beta-1}$, $b_{n}=\left[  n\right]
_{q}^{\beta}$, $0<\beta\leq\frac{2}{3}$, $n\in\mathbb{N}$ and $x\neq-\frac
{1}{a_{n}}$.
\end{definition}

In the case $q=1$ these polynomials coincide with the classical ones. For
$q\neq1$ one gets a new class of polynomials having interesting properties. It
should be mentioned that in the case $q\in\left(  0,1\right]  $ $q$%
-Bal\'{a}zs-Szabados operators generate positive linear operators
$R_{n,q}:f\rightarrow R_{n,q}\left(  f;x\right)  $. In the case $q>1$
positivity fails, however, the results of this paper show that in this case
approximating properties of $q$-Bal\'{a}zs-Szabados operators may be better
than the case $0<q\leq1$.

Throughout this paper, let $\mathbb{D}_{R}:=\left\{  z\in\mathbb{C}:\left\vert
z\right\vert <R\right\}  $ denote the disk of radius $R$ centered at $0.$
Moreover, it is assumed that $f\left(  z\right)  =%
{\displaystyle\sum\nolimits_{m=0}^{\infty}}
c_{m}z^{m},\ $for $z\in\mathbb{D}_{R}$.

Assuming $f:\mathbb{D}_{R}\cup\left[  R,+\infty\right)  \rightarrow\mathbb{C}$
and simply replacing $x$ by $z$ in (\ref{bs}) we obtain the complex form of
the $q$-Bal\'{a}zs-Szabados operator:
\[
R_{n,q}\left(  f;z\right)  =\dfrac{1}{\left(  1+a_{n}z\right)  ^{n}}\sum
_{k=0}^{n}f\left(  \dfrac{\left[  k\right]  _{q}}{b_{n}}\right)  \left[
\begin{array}
[c]{c}%
n\\
k
\end{array}
\right]  _{q}\left(  a_{n}z\right)  ^{k}%
{\displaystyle\prod\limits_{s=0}^{n-k-1}}
\left(  1+\left(  1-q\right)  \left[  s\right]  _{q}a_{n}z\right)  ,
\]
where again $a_{n}=\left[  n\right]  _{q}^{\beta-1}$, $b_{n}=\left[  n\right]
_{q}^{\beta}$, $0<\beta\leq\frac{2}{3}$, $n\in\mathbb{N}$, $z\in\mathbb{C}$
and $z\neq-\frac{1}{a_{n}}$.

\begin{remark}
The complex operators $R_{n,q}\left(  f;z\right)  $ are well defined and
analytic for all $n\geq n_{0}$ and $\left\vert z\right\vert \leq r<\left[
n_{0}\right]  _{q}^{1-\beta}$. Indeed, in this case we easily obtain that
$z\neq-\frac{1}{a_{n}}$, for all $\left\vert z\right\vert \leq r<\left[
n_{0}\right]  _{q}^{1-\beta}$ and $n\geq n_{0}$, which implies that $1/\left(
1+a_{n}z\right)  ^{n}$ is analytic.
\end{remark}

\begin{remark}
There exists a close connection between $R_{n,q}\left(  f;z\right)  $ and the
complex $q$-Bernstein polynomials given by%
\[
B_{n,q}\left(  f;z\right)  =\sum_{k=0}^{n}f\left(  \dfrac{\left[  k\right]
_{q}}{\left[  n\right]  _{q}}\right)  \left[
\begin{array}
[c]{c}%
n\\
k
\end{array}
\right]  _{q}z^{k}%
{\displaystyle\prod\limits_{s=0}^{n-k-1}}
\left(  1-q^{s}z\right)  .
\]
Indeed, denoting $F_{n}\left(  z\right)  =f\left(  \dfrac{\left[  n\right]
_{q}}{b_{n}}z\right)  $, we easily get%
\[
R_{n,q}\left(  f;z\right)  =B_{n,q}\left(  F_{n};\dfrac{a_{n}z}{1+a_{n}%
z}\right)  ,
\]
valid for all $n\geq n_{0}$ and $\left\vert z\right\vert \leq r<\left[
n_{0}\right]  _{q}^{1-\beta}.$ This connection will be essential in our
reasonings. For monomials $f\left(  z\right)  =e_{m}\left(  z\right)  =z^{m}$
it can be written as follows:%
\[
R_{n,q}\left(  e_{m};z\right)  =\left[  n\right]  _{q}^{m\left(
1-\beta\right)  }B_{n,q}\left(  e_{m};\dfrac{a_{n}z}{1+a_{n}z}\right)  .
\]

\end{remark}

\begin{remark}
The lack of positivity makes the investigation of convergence in the case
$q>1$ essentially more difficult than for $0<q<1$. Notice that, the complex
$q$-Bernstein type operators in the case $q>1$ systematically are studied in
\cite{ostrov1}, \cite{ostrov2}, \cite{wu1}, \cite{wu2}, \cite{mah1}, and
\cite{mah2}.
\end{remark}

\begin{remark}
Approximation properties of the complex Bal\'{a}zs-Szabados operators are
studied in \cite{galbook1}. Notice that unlike to \cite{galbook1} the growth
conditions of exponential-type on $f$ is omitted. The only condition imposed
to $f$ is to be uniformly continuous and bounded on $\left[  0,+\infty\right)
.$ Therefore our results are new even for the classical Bal\'{a}zs-Szabados operators.
\end{remark}

\begin{theorem}
\label{upperestimate}Let $n_{0}\geq2$, $0<\beta\leq\dfrac{2}{3}$. Assume that
$f:\mathbb{D}_{R}\cup\left[  R,+\infty\right)  \rightarrow\mathbb{C}$ is
uniformly continuous and bounded on $\left[  0,+\infty\right)  $, is analytic
in $\mathbb{D}_{R}.$ Then%
\begin{align*}
\left\vert R_{n,q}\left(  f;z\right)  -f\left(  z\right)  \right\vert  &
\leq\frac{1}{\left[  n\right]  _{q}^{\beta}}%
{\displaystyle\sum\nolimits_{m=2}^{\infty}}
\left\vert c_{m}\right\vert m\left(  m-1\right)  \left(  4q^{2}r\right)
^{m}+\frac{2r}{\left[  n\right]  _{q}^{1-\beta}}%
{\displaystyle\sum\nolimits_{m=1}^{\infty}}
\left\vert c_{m}\right\vert \left(  2r\right)  ^{m},\ \\
q  &  \geq1,\ \ \frac{1}{2}<r<\frac{R}{4q^{2}}\leq\frac{1}{2}\left[
n_{0}\right]  _{q}^{1-\beta}.
\end{align*}
\bigskip
\end{theorem}

In \cite{totik}, Totik settled the saturation properties of $R_{n}$. Among
other things he proved the Voronovskaja-type result for $0<\beta\leq\frac
{1}{2},$ $\beta\geq\frac{2}{3}$. The complete asymptotic expansion for $R_{n}$
is given by Abel and Della Vecchia \cite{abel}.

Next, we study Voronovskaja type formulas of the $q$-Bal\'{a}zs-Szabados
operators of a function $f$ analytic in the disc $\mathbb{D}_{R}$. In order to
formulate Voronovskaja type theorem we define the following function%
\begin{equation}
L_{q}^{\beta}\left(  f;z\right)  :=\left\{
\begin{tabular}
[c]{lll}%
$\dfrac{D_{q}f\left(  z\right)  -f^{\prime}\left(  z\right)  }{q-1},$ & if &
$\left\vert z\right\vert <R/q,\ R>q>1,\ \ 0<\beta<\frac{1}{2},$\\
$-z^{2}f^{\prime}\left(  z\right)  +\dfrac{D_{q}f\left(  z\right)  -f^{\prime
}\left(  z\right)  }{q-1},$ & if & $\left\vert z\right\vert
<R/q,\ R>q>1,\ \ \beta=\frac{1}{2},$\\
$-z^{2}f^{\prime}\left(  z\right)  ,$ & if & $\left\vert z\right\vert
<R/q,\ R>q>1,\ \ \ \frac{1}{2}<\beta<1,$%
\end{tabular}
\ \ \ \right.  \label{lq}%
\end{equation}
\bigskip and for $q=1,$%
\[
L_{q}^{\beta}\left(  f;z\right)  :=\left\{
\begin{tabular}
[c]{lll}%
$\dfrac{z}{2}f^{\prime\prime}\left(  z\right)  ,$ & if & $\left\vert
z\right\vert <R,\ \ 0<\beta<\frac{1}{2},$\\
$-z^{2}f^{\prime}\left(  z\right)  +\dfrac{z}{2}f^{\prime\prime}\left(
z\right)  ,$ & if & $\left\vert z\right\vert <R,\ \ \beta=\frac{1}{2},$\\
$-z^{2}f^{\prime}\left(  z\right)  ,$ & if & $\left\vert z\right\vert
<R,\ \ \ \frac{1}{2}<\beta<1.$%
\end{tabular}
\ \ \ \right.
\]

In the case of complex variable, the qualitativeVoronovskaja-type result for
$R_{n}$ is proved by Gal \cite{galbook1}. Note that the case $\beta=\frac
{1}{2}$ remained open in \cite{galbook1}. We prove the following quantitative
Voronovskaja type theorem for $R_{n,q},$ which covers the case $\beta=\frac
{1}{2}$. Moreover, our results are new for the classical Bal\'{a}zs-Szabados
operators ($q=1$).\bigskip

\begin{theorem}
\label{vor}Let $n_{0}\geq2$, $0<\beta\leq\dfrac{2}{3}$. Assume that
$f:\mathbb{D}_{R}\cup\left[  R,+\infty\right)  \rightarrow\mathbb{C}$ is
uniformly continuous and bounded on $\left[  0,+\infty\right)  $, is analytic
in $\mathbb{D}_{R}.$ Then

\begin{enumerate}
\item[(i)] For $0<\beta<\frac{1}{2},$ $\frac{1}{2}<r<\dfrac{R}{\max\left(
4q,2q^{2}\right)  }<R\leq\frac{1}{2}\left[  n_{0}\right]  _{q}^{1-\beta},$
$\left\vert z\right\vert \leq r$, we have
\begin{align*}
&  \left\vert R_{n,q}\left(  f;z\right)  -f\left(  z\right)  -\dfrac
{1}{\left[  n\right]  _{q}^{\beta}}\dfrac{z\left(  D_{q}f\left(  z\right)
-f^{\prime}\left(  z\right)  \right)  }{q-1}\right\vert \\
&  \leq\frac{4}{\left[  n\right]  _{q}^{2\beta}}\sum_{m=0}^{\infty}\left\vert
c_{m}\right\vert \left(  m-2\right)  (4qr)^{m-2}+\dfrac{4}{\left[  n\right]
_{q}^{1-\beta}}\sum_{m=0}^{\infty}\left\vert c_{m}\right\vert m\left(
m-1\right)  (2q^{2}r)^{m+1};
\end{align*}

\item[(ii)] For $\frac{1}{2}<\beta\leq\dfrac{2}{3},$ $\frac{1}{2}<r<\dfrac
{R}{4q}<R\leq\frac{1}{2}\left[  n_{0}\right]  _{q}^{1-\beta},$ $\left\vert
z\right\vert \leq r$, we have%
\[
\left\vert R_{n,q}\left(  f;z\right)  -f\left(  z\right)  +\dfrac{1}{\left[
n\right]  _{q}^{1-\beta}}z^{2}f^{\prime}\left(  z\right)  \right\vert
\leq\dfrac{6}{\left[  n\right]  _{q}^{\beta}}\sum_{m=0}^{\infty}\left\vert
c_{m}\right\vert m\left(  m-1\right)  (4qr)^{m};
\]

\item[(iii)] For $\beta=\frac{1}{2},$ $\frac{1}{2}<r<\dfrac{R}{4q^{2}}%
<R\leq\frac{1}{2}\left[  n_{0}\right]  _{q}^{1-\beta},$ $\left\vert
z\right\vert \leq r$, we have%
\begin{align*}
&  \left\vert R_{n,q}\left(  f;z\right)  -f\left(  z\right)  +\dfrac{1}%
{\sqrt{\left[  n\right]  _{q}}}z^{2}f^{\prime}\left(  z\right)  -\dfrac
{1}{\sqrt{\left[  n\right]  _{q}}}\dfrac{z\left(  D_{q}f\left(  z\right)
-f^{\prime}\left(  z\right)  \right)  }{q-1}\right\vert \\
&  \leq\dfrac{9}{\left[  n\right]  _{q}}\sum_{m=0}^{\infty}\left\vert
c_{m}\right\vert m^{2}\left(  m-1\right)  ^{2}\left(  4q^{2}r\right)  ^{m}.
\end{align*}

\end{enumerate}
\end{theorem}

By using the above Voronovskaja's theorem, we will obtain the exact order in
approximation by the complex $q$-Bal\'{a}zs-Szabados operators. In this sense,
we present the following results.

\begin{theorem}
\label{degree}Let $n_{0}\geq2$, $0<\beta\leq\dfrac{2}{3}$. Assume that
$f:\mathbb{D}_{R}\cup\left[  R,+\infty\right)  \rightarrow\mathbb{C}$ is
uniformly continuous and bounded on $\left[  0,+\infty\right)  $, is analytic
in $\mathbb{D}_{R}$.

\begin{enumerate}
\item[(i)] If $0<\beta<\frac{1}{2},$ $\frac{1}{2}<r<\dfrac{R}{\max\left(
4q,2q^{2}\right)  }<R\leq\frac{1}{2}\left[  n_{0}\right]  _{q}^{1-\beta},$ and
$f$ is not a polynomial of degree $\leq1$ in $\mathbb{D}_{R}$, then
\[
\left\Vert R_{n,q}\left(  f\right)  -f\right\Vert _{r}\sim\frac{1}{\left[
n\right]  _{q}^{\beta}},\ \ n\in\mathbb{N}.
\]

\item[(ii)] If $\frac{1}{2}<\beta\leq\dfrac{2}{3},$ $\frac{1}{2}<r<\dfrac
{R}{4q}<R\leq\frac{1}{2}\left[  n_{0}\right]  _{q}^{1-\beta},$ and $f$ is not
a constant function in $\mathbb{D}_{R}$, then
\[
\left\Vert R_{n,q}\left(  f\right)  -f\right\Vert _{r}\sim\frac{1}{\left[
n\right]  _{q}^{1-\beta}},\ \ n\in\mathbb{N}.
\]

\item[(iii)] For $\beta=\frac{1}{2},$ $\frac{1}{2}<r<\dfrac{R}{4q^{2}}%
<R\leq\frac{1}{2}\left[  n_{0}\right]  _{q}^{1-\beta},$ and $f$ is not a
constant function in $\mathbb{D}_{R}$, then
\[
\left\Vert R_{n,q}\left(  f\right)  -f\right\Vert _{r}\sim\frac{1}%
{\sqrt{\left[  n\right]  _{q}}},\ \ n\in\mathbb{N}.
\]

\end{enumerate}
\end{theorem}


\bigskip

\begin{thebibliography}{99}                                                                                               %


\bibitem {balazs1}K. Bal\'{a}zs, Approximation by Bernstein type rational
functions, Acta Math. Acad. Sci. Hungar., 26, (1975) 123-134.

\bibitem {balazs2}K. Bal\'{a}zs and J. Szabados, Approximation by Bernstein
type rational functions, II, Acta Math. Acad. Sci. Hungar., 40 (1982) 3-4, 331-337.

\bibitem {abel}U. Abel and B. Della Vecchia,. Asymptotic approximation by the
operators of K. Bal azs and Szabados, Acta Sci. Math.(Szeged), 66, (2000) No.
1-2, 137\{145.

\bibitem {dogru}O. Dogru, On statistical approximation properties of Stancu
type bivariate generalization of $q$-Bal\'{a}zs-Szabados operators, Proc. of
Int. Conf. on Numer. Anal. and Approx. Th. Cluj-Napoca, Romanya, 2002.

\bibitem {galbook1}S.G. Gal, Approximation by Complex Bernstein and
Convolution Type Operators. World Scientific Publishing, New Jersey, London,
Singapore, Beijing, Shanghai, Hong Kong, Taipei, Chennai (2009)

\bibitem {ozkan}N. Ispir and Y. \"{O}zkan, Approximation properties of complex
$q$-Bal\'{a}zs-Szabados operators in compact disks, Journal of Inequalities
and Applications 2013, 2013:361

\bibitem {totik}V. Totik, Saturation for Bernstein-type rational functions,
Acta Math. Hungar., 43 (1984), 219--250.

\bibitem {agarwal}R. P Agarwal\ and V. Gupta, On $q$-analogue of a complex
summation-integral type operators in compact disks, Journal of Inequalities
and Applications 2012, 2012:111.

\bibitem {AAR}Andrews G E, Askey R, Roy R. Special functions. Cambridge:
Cambridge University Press; 1999.

\bibitem {ostrov1}S. Ostrovska: $q$-Bernstein polynomials and their iterates.
J. Approximation Theory 123 (2003), 232--255.

\bibitem {ostrov2}S. Ostrovska: The sharpness of convergence results for
$q$-Bernstein polynomials in the case $q>1$. Czech. Math. J. 58 (2008), 1195--1206.

\bibitem {wu1}H. Wang, X. Wu: Saturation of convergence for $q$-Bernstein
polynomials in the case $q>1$. J. Math. Anal. Appl. 337 (2008), 744--750.

\bibitem {wu2}Z. Wu: The saturation of convergence on the interval $[0,1]$ for
the $q$-Bernstein polynomials in the case $q>1$. J. Math. Anal. Appl. 357
(2009), 137--141.

\bibitem {mah1}N. I. Mahmudov, Approximation by $q$-Durrmeyer type polynomials
in compact disks in the case $q>1$. Appl. Math. Comput. 237 (2014), 293--303

\bibitem {mah2}N. I. Mahmudov, Approximation by genuine $q$%
-Bernstein-Durrmeyer polynomials in compact disks in the case $q>1$. Abstr.
Appl. Anal. 2014, Art. ID 959586, 11 pp.
\end{thebibliography}
\end{document}